\documentclass{amsart}
\usepackage{amssymb}
\usepackage{amscd}

\usepackage[
hyperref
]{degt}

\makeatletter

\def\Prod{\BG3\cdot\Inn\AA}

\let\rdeg\width  
\let\bdeg\dg     

\def\bGb{\bar\Gb}

\def\scirc{^\sharp}

\def\be{\bold e}
\def\bu{\bold u}
\def\bv{\bold v}

\def\CR{\Cal R}

\def\CV{\Cal V}
\def\CU{\Cal U}

\let\BB=B  
\let\CC=C  
\let\EE=E  
\let\DD=D  
\let\disk\Omega
\let\slope\varkappa

\def\BM{\frak{Im}} 
\let\AM=A          
\def\AM{\mathrm{A}} 
\pdef\bAM{\bar\AM}
\def\VM{\Cal V}
\pdef\bVM{\bar\VM}

\let\pp\pi               

\def\FF#1{\Bbb F_{#1}}   
\def\AA{\Bbb G}          
\def\AA{\frak F}         
\let\cp\Phi              
\def\bk{\bar\Bbbk}

\def\Espec{\bold S}      
\def\Rtriv{\bold R}      
\def\CE{\Cal E}

\def\Supp{\QOPNAME{Supp}}

\def\bm{\frak m}
\let\kk\kappa

\def\inserthyphen{\ifcat\next a-\fi\ignorespaces}
\let\BLACK\bullet
\let\WHITE\circ
\def\CROSS{\vcenter{\hbox{$\scriptstyle\mathord\times$}}}
\let\STAR*
\def\TRIANG{\vcenter{\hbox{$\scriptstyle\mathord\vartriangle$}}}
\def\pblack-{$\BLACK$\futurelet\next\inserthyphen}
\def\pwhite-{$\WHITE$\futurelet\next\inserthyphen}
\def\pcross-{$\CROSS$\futurelet\next\inserthyphen}
\def\pstar-{$\STAR$\futurelet\next\inserthyphen}
\def\ptriang-{$\TRIANG$\futurelet\next\inserthyphen}
\def\black{\protect\pblack}
\def\white{\protect\pwhite}

\def\NO#1{n_{#1}}
\def\nblack{\NO\BLACK}
\def\nwhite{\NO\WHITE}

\let\Bs\beta
\def\Bs{\bar\Gs}
\def\Bb{\bar\Gb}
\def\bG{\bar G}

\let\mp\psi  

\def\SkGamma{\BLACK{\joinrel\relbar\joinrel\relbar\joinrel}\WHITE}

\let\vI\I
\let\vII\II
\let\vIII\III
\let\vIV\IV

\let\vZ0

\def\ee{\qopname@{e}}

\def\reg(#1){(\!(#1)\!)}

\def\Maple{{\tt Maple}}

\newcounter{line}

\makeatother

\title{The Alexander module of a trigonal curve. II}

\author{Alex Degtyarev}

\address{%
Bilkent University\\
Department of Mathematics\\
06800 Ankara, Turkey}

\email{degt@fen.bilkent.edu.tr}



\begin{document}

\begin{abstract}
We complete the
enumeration of the possible roots of the Alexander polynomial (both
conventional and over finite fields) of a trigonal curve.
The curves are not assumed proper or irreducible.
\end{abstract}

\maketitle

\section{Introduction}

This paper is a continuation of~\cite{degt:Alexander}: we complete the
enumeration of the possible roots of the Alexander polynomial (both
conventional and over finite fields) of a trigonal curve.
Unlike~\cite{degt:Alexander}, here we do not assume the curves irreducible,
as this assumption does not lead to an essential simplification of the
results.

An emphasis is given to the \emph{exceptional} roots~$\xi$,
\ie, those with the
multiplicative order $N:=\ord(-\xi)\ge7$. Such roots are \emph{not}
controlled by congruence subgroups of the modular group.

Since the paper is a sequel, we only recall very briefly the necessary
notions and preliminary results, concentrating on
the explanation of the new
approach that lets us improve the estimates found in~\cite{degt:Alexander}.
For all details, further speculations, and references
concerning the history of the subject
and the previously known results on the
Alexander module/polynomial of an algebraic
curve in an algebraic surface, the reader is directed
to~\cite{degt:Alexander} and~\cite{degt:trigonal}.

\subsection{Principal results}
Let $\CC\subset\Sigma_d$ be a trigonal curve in a Hirzebruch surface,
see \autoref{s.tc} for the definitions,
and consider the affine fundamental group
\[*
\piaff{\CC}:=\pi_1(\Sigma_d\sminus(\CC\cup\EE\cup F_\infty)),
\]
where $F_\infty$ is a generic fiber of~$\Sigma_d$. and $E$ is the
exceptional section.
There is a natural epimorphism $\deg\:\piaff{\CC}\onto\Z$ sending a meridian of
a tubular neighborhood of~$\CC$ to $1\in\Z$.
The abelianization $\AM_\CC$ of the kernel $\Ker{\deg}$ is called the
\emph{Alexander module} (or \emph{Alexander invariant}) of~$\CC$.
This group is indeed a module over the ring $\Lambda:=\Z[t,t\1]$ of Laurent
polynomials; the action of~$t$ is given by
$[h]\mapsto[aha\1]$, where $h\in\Ker\deg$ and $a\in\piaff{\CC}$ is any element
of degree~$1$.
Alternatively, $\AM_\CC=H_1(X)$ is the homology of the infinite cyclic covering
$X\to\Sigma_d\sminus(\CC\cup\EE\cup F_\infty)$ corresponding to~$\deg$,
and the action
of~$t$ is induced by the deck translation of the covering.

Denote $\Bbbk_0:=\Q$ and $\Bbbk_p:=\FF{p}$ (the field with $p$ elements) if
$p$ is a prime. Unless $\CC$ is isotrivial,
the product $\AM_\CC\otimes\Bbbk_p$ is a torsion module
over the principal ideal domain $\Lambda\otimes\Bbbk_p$;
its order $\Delta_{\CC,p}\in\Lambda\otimes\Bbbk_p$
is called the \emph{$({\bmod}\,p)$-Alexander polynomial} of~$\CC$.
We are interested in the roots of $\Delta_{\CC,p}$.
More precisely,
let~$\xi$ be an algebraic number over~$\Bbbk_p$, denote by
$\mp_\xi\in\Lambda\otimes\Bbbk_p$ its minimal polynomial, and consider the
minimal field
$\bk:=\Lambda(\xi):=(\Lambda\otimes\Bbbk_p)/\mp_\xi$ containing~$\xi$.
Then $\AM_\CC(\xi):=(\AM_\CC\otimes\Bbbk_p)/\mp_\xi$ is a $\bk$-vector
space, and we are interested in the pairs $(p,\xi)$ for which this space may
have positive dimension, \ie, $\mp_\xi$ may appear as a factor of the
Alexander polynomial $\Delta_{\CC,p}$ of a non-isotrivial trigonal curve.
According to~\cite{degt:Alexander}, the multiplicative order
$N:=\ord(-\xi)$ must be finite, and the principal result of the present paper
is the following theorem,
which is proved in \autoref{S.proof.main}, see \autoref{s.main}.

\table
\caption{Exceptional factors of~$\Delta$, $N\ge7$%
 \noaux{ (see \autoref{rem.factors})}}\label{tab.factors}
\bgroup
\def\*{\llap{$^{*}$}}%
\let\comma,\catcode`\,\active\def,{$\comma\ \ $}%
\let\scolon;\catcode`\;\active\def;{$\scolon\ \ $}%
\def\get[#1]{$#1$\hss}%
\def\getG(#1;#2;#3;#4){$(#1\scolon#2\comma#3\scolon#4)$\hss}%
\centerline{\vbox{\halign{\strut\hss\refstepcounter{line}#\theline&
 \quad\hss$#$&\quad\hss$#$&\quad\get#&
 \quad\hss\getG#\hss\cr
\omit\strut\hss$\#$&\omit\hss$p$&$N$&\omit\qquad
 Factors $\mp_\xi\in\FF{p}[t]$ of $\Delta_{C\comma p}$\hss&\omit\quad\hss$\bG\subset\MG$\hss\cr
\noalign{\vskip2pt \hrule\vskip3pt}
\*&  2& 7&[t^3+t+1, t^3+t^2+1]&(9;1;0;1^2 7^1)\cr
\*&   &15&[t^4+t+1, t^4+t^3+1]&(17;1;2;1^2 15^1)\cr
\*&  3&8&[t^2+2t+2, t^2+t+2]&(10;0;1;1^2 8^1)\cr
\*&  5&8&[t^2+2, t^2+3]&(78;0;0;1^6 8^9)\cr
&   &  12&[t^2+2t+4, t^2+3t+4]&(52;0;4;1^4 12^4)\cr
\*& 11&10&[t+2; t+6; t+7; t+8]&(24;2;0;1^2 2^1 10^2)\cr
\*& 13&12&[t+2, t+7; t+6, t+11]&(14;0;2;1^2 12^1)\cr
\*& 17&8&[t+2, t+9; t+8, t+15]&(36;0;0;1^4 8^4)\cr
& 19&   9&[t+4, t+5; t+6, t+16; t+9, t+17]&(20;0;2;1^2 9^2)\cr
&   &  18&[t+2; t+3; t+10; t+13; t+14; t+15]&(40;2;4;1^2 2^1 18^2)\cr
\*& 29&7&[t+7, t+25; t+16, t+20; t+23, t+24]&(60;0;0;1^4 7^8)\cr
& 37&   9&[t+7, t+16; t+9, t+33; t+12, t+34]&(76;0;4;1^4 9^8)\cr
\*& 43&7&[t+4, t+11; t+16, t+35; t+21, t+41]&(132;0;0;1^6 7^{18})\cr
\crcr}}}\egroup
\endtable

\theorem\label{th.main}
Let $p$ be a prime or zero, and
assume that the $({\bmod}\,p)$-Alexander polynomial
of a
non-isotrivial \rom(see \autoref{s.tc}\rom)
trigonal curve~$\CC$ has a root $\xi\in\bk\supset\Bbbk_p$
such that $N:=\ord(-\xi)\ge7$. Then the pair $(p,\mp_\xi)$ is one of those
listed in \autoref{tab.factors}.
The pairs $(p,\mp_\xi)$ marked with a $^*$ in the table
do appear in the Alexander polynomials of
\emph{proper} trigonal curves\rom; the others do not.
\endtheorem


\remark[comments to \autoref{tab.factors}]\label{rem.factors}
Listed in the table are triples $(p,N,\mp_\xi)$ and, for each triple,
certain information about the skeleton (see~\autoref{s.skeletons})
of the corresponding universal
subgroup (see \autoref{s.universal}):
$(e;v_\WHITE,v_\BLACK;r)$, where $e$ is the number of edges, $v_\WHITE$
and~$v_\BLACK$ are the numbers of monovalent \white-- and \black-vertices,
respectively, and $r$ is the set of the region widths in the partition
notation. Note that these data do not determine the skeleton uniquely: in
fact, the polynomials~$\mp_\xi$ with isomorphic skeletons are separated by
commas rather than semicolons in the table.
\endremark

\remark
The case $N>10$ is settled in~\cite{degt:Alexander}, where the rest of
\autoref{th.main} (the range $7\le N\le10$)
was conjectured.
All triples $(p,N,\mp_\xi)$ with $N\le6$ are also
enumerated in~\cite{degt:Alexander}, see \autoref{s.N<=6} for
a few further details.
\endremark

\addendum\label{ad.<=1}
In all cases listed in \autoref{tab.factors}, the module $\AM_\CC(\xi)$ has
a geometric presentation of the form $\AM(\xi)/\bk\be_1$,
see \autoref{s.Alexander} and \autoref{s.Burau} for the definition
of~$\AM$ and~$\be_i$\rom;
in particular, $\dim_{\bk}\AM_\CC(\xi)=1$.
Furthermore, at most one factor as in the table may appear in the Alexander
polynomial of any particular curve.
\endaddendum

This statement is proved in \autoref{proof.<=1}.

\subsection{Open problems}
The approach chosen here and in~\cite{degt:Alexander}, namely, the study of
the specializations $\AM_C(\xi)$, simplifies the original question about the
structure of the $\Lambda$-module $\AM_C$.
In particular, we ignore the
higher torsion of the form
$\CG{p^r}$ and $(\Lambda\otimes\FF{p})/\mp^r$, $p>0$.
In general, it is clear that $\AM_C$ is annihilated by a polynomial of the
form $t^n-1$, and it appears reasonable to study the quotients
$\AM_C/\cp_N(-t)$, where $\cp_N$ is the cyclotomic polynomial of order~$N$.
(Still, since $\Lambda$ is not a principal ideal domain,
extra work is needed to recover~$\AM_C$ from these quotients.)
At present, the structure of the $\Lambda/\cp_N(-t)$-module
$\AM_C/\cp_N(-t)$ is known for $N=1$,
see~\cite{degt:trigonal}, and for $2\le N\le6$, see \autoref{s.N<=6}.
If $N\ge7$, using \autoref{th.main} and \autoref{ad.<=1}, one can only state
that $\AM_C/\cp_N(-t)$ is generated over~$\Lambda$ by a single element and
that it is a finite abelian $p$-group for some value
of prime~$p$ (compatible with that of~$N$) as in \autoref{tab.factors}.

Examples of modules $\AM_C$ with higher torsion
$\CG{p^r}$ or $(\Lambda\otimes\FF{p})/\mp^r$ are known,
see~\cite{degt:trigonal} and \autoref{s.N<=6}. However, in all these examples
one has $N\le6$; the case $N\ge7$ has never been studied from this point of
view.

Another open question is the realizability of all modules discovered by
trigonal curves. Our proof of \autoref{th.main} is in the group-theoretical
settings, and each pair $(p,\mp_\xi)$ listed in \autoref{tab.factors}
\emph{is} realized by a certain genus zero subgroup~$G$ of the Burau
group~$\Bu3$. If $G\subset\BG3$, it is in turn realized
as the monodromy group of a proper trigonal
curve due to \autoref{th.groups}; the other cases, leading to improper
curves, are to be the subject of a further investigation.

\subsection{Contents of the paper}
The first two sections contain preliminary material: we site a few notions
and results needed in the proof of the main theorem.
In \autoref{S.groups}, we introduce the braid group~$\BG3$ and its extension
$\Prod$ and
recall the reduced Burau representation (see \autoref{s.Burau}) and
its further specialization to the modular group $\MG:=\PSL(2,\Z)$.
Crucial for the sequel is the description of subgroups of~$\MG$ in terms of
\emph{skeletons} (certain bipartite ribbon graphs, see \autoref{s.skeletons})
and the description of their lifts to $\BG3$ or~$\Bu3$ by means of the
\emph{type specification} (see \autoref{s.type}).
In \autoref{S.tc}, we discuss trigonal curves in Hirzebruch surfaces
(see \autoref{s.tc}),
the braid monodromy and the monodromy group of such a curve
(see \autoref{s.monodromy}),
and the Zariski--van Kampen theorem (see \autoref{s.vanKampen}).
Then we introduce the conventional and extended Alexander module of a
trigonal curve, expressing them in terms of the monodromy group
(see \autoref{s.Alexander})
and introduce the universal subgroup corresponding
to a given Alexander module (see \autoref{s.universal}).
For completeness, we extend to all,
not necessarily proper or irreducible, trigonal
curves the results of~\cite{degt:Alexander} concerning the roots~$\xi$ of the
Alexander polynomial with $N:=\ord(-\xi)\le6$ (see \autoref{s.N<=6}).

\autoref{th.main} is proved in \autoref{S.proof.main}.
We recall the computation of the local modules and the estimate $N\le26$
found in~\cite{degt:Alexander} (see \autoref{s.local}), show that
$\dim_{\bk}\AM_C(\xi)\le1$ whenever $N\ge7$ (see \autoref{s.VM=0}), and
engage into improving the
above estimate to $N\le6$ with the exception of finitely
many explicitly listed cases (see \autoref{s.reduction}): the strategy
consists in showing, by computing a number of resultants, that the assumption
$N\ge7$ implies that a certain ribbon graph has too many too large regions
and thus cannot be planar.
The exceptional cases are eliminated by the explicit computation of the
corresponding universal groups, which are all finite (see
\autoref{s.groups}), and this concludes the proof of \autoref{th.groups},
which restates \autoref{th.main} in group-theoretical terms.
The computation in \autoref{s.reduction} and \autoref{s.groups} is heavily
computer aided;
it was done using \Maple,
and we only outline the approach.
The formal reduction of \autoref{th.main} to \autoref{th.groups}
is explained in \autoref{s.main},
and \autoref{ad.<=1} is proved in \autoref{proof.<=1}.

As usual, ends of proofs are marked with \qedsymbol. Statements whose proofs
are omitted are marked with either \donesymbol\! or \pnisymbol\!. In the
former case, the proof is either trivial or already explained; in the latter
case, the reader is directed to the literature, which is usually referred to
in the header of the statement.

\subsection{Acknowledgements}
I would like to thank the organizers and participants of the Fukuoka
Symposium for the great job that resulted in a remarkable meeting.
This paper was
essentially
written during my sabbatical stay at
\emph{l'Universit\'{e} Pierre et Marie Curie} (Paris 6); I am grateful to this
institution for its hospitality. I would also like to extend my gratitude
to Mouadh Akriche for his helpful comments.

\section{The groups $\BG3$, $\Bu3$, and \pdfstr{Gamma}{$\MG$}}\label{S.groups}

In this section, se recall briefly a few necessary fact concerning the braid
group~$\BG3$, the modular group $\MG:=\PSL(2,\Z)$, and their relation to
bipartite ribbon graphs.
A detailed treatment of the latter subject is found in~\cite{degt:book};
here, we merely recall the basic definitions and explain briefly the
geometric insight. For all proofs, the reader is also referred
to~\cite{degt:book}.

\subsection{The braid group $\BG3$}\label{s.BG3}

Artin's \emph{braid group}~$\BG3$ on three strands is the group
\[*
\BG3:=\<\Gs_1,\Gs_2\,|\,\Gs_1\Gs_2\Gs_1=\Gs_2\Gs_1\Gs_2\>.
\]
The generators $\Gs_1,\Gs_2$ above are called \emph{Artin generators}.
There is an epimorphism $\bdeg\:\BG3\onto\Z$, $\Gs_1,\Gs_2\mapsto1\in\Z$,
called the \emph{degree}. Furthermore, there is a canonical
faithful representation
$\BG3\to\Aut\<\Ga_1,\Ga_2,\Ga_3\>$, the Artin generators acting \via
\[
\Gs_1\:\Ga_1\mapsto\Ga_1\Ga_2\Ga_1\1,\quad \Ga_2\mapsto\Ga_1;\qquad
\Gs_2\:\Ga_2\mapsto\Ga_2\Ga_3\Ga_2\1,\quad \Ga_3\mapsto\Ga_2.
\label{eq.braid.action}
\]
According to Artin~\cite{Artin}, $\BG3$ can be identified with the subgroup
of $\Aut\<\Ga_1,\Ga_2,\Ga_3\>$ consisting of the automorphisms taking each
generator to a conjugate of a generator and preserving the product
$\Gr:=\Ga_1\Ga_2\Ga_3$.

Throughout the paper, we reserve the notation $\AA$ for the free group
$\<\Ga_1,\Ga_2,\Ga_3\>$ equipped with a distinguished $\BG3$-orbit of bases,
called \emph{geometric}.
(Note that each particular geometric basis gives rise to its own pair
$\Gs_1,\Gs_2$ of
Artin generators of~$\BG3$.)
The element
$\Gr:=\Ga_1\Ga_2\Ga_3\in\AA$ does not depend on the choice of a geometric
basis. The center of~$\BG3$ is the infinite cyclic group generated
by~$\Delta^2$, where $\Delta:=\Gs_1\Gs_2\Gs_1=\Gs_2\Gs_1\Gs_2$ is the
\emph{Garside element}; one has $\Delta^2(\Ga)=\Gr\Ga\Gr\1$ for any
$\Ga\in\AA$.
The Garside element~$\Delta$ depends on the choice of a geometric basis,
whereas $\Delta^2$ does not.

There is a well defined \emph{degree epimorphism}
$\deg\:\AA\onto\Z$ taking each geometric generator to $1\in\Z$; it is
preserved by~$\BG3$.

With improper trigonal curves in mind, we will consider a larger
group~$\Prod$, where $\Inn\AA\subset\Aut\AA$ is the subgroup of inner
automorphisms. Since $\Inn\AA$ is normal in $\Aut\AA$, the product is indeed
a subgroup. One has $\BG3\cap\Inn\AA=\<\Delta^2\>$; hence,
the map $\bdeg\:\BG3\onto\Z$ extends to $\Prod$ \via\
$\bdeg(\Gb\Ga)=\bdeg\Gb+2\deg\Ga$, where $\Gb\in\BG3$
and $\Ga\in\Inn\AA\cong\AA$.
(We consider the left adjoint action $\Ga(\Ga')=\Ga\Ga'\Ga\1$; observe that
$\bdeg\Delta^2=6=2\deg\Gr$.)

\subsection{The Burau representation}\label{s.Burau}

The \emph{universal Alexander module} is the abelianization~$\AM$ of the
kernel $\Ker[\deg\:\AA\onto\Z]$. It is a module over the ring
$\Lambda:=\Z[t,t\1]$ of integral Laurent polynomials, $t$ acting \via\
$t[h]=[\Ga h\Ga\1]$, where $h\in\Ker\deg$ and $\Ga\in\AA$ is any element of
degree~$1$.
A simple computation using the Reidemeister--Schreier method shows that
$\AM=\Lambda\be_1\oplus\Lambda\be_2$, where
$\be_i:=[\Ga_{i+1}\Ga_i\1]$, $i=1,2$;
these generators depend on the choice of a geometric basis.

Since the action of $\Prod$ preserves the degree, it descends to an action
on~$\AM$, giving rise to a representation $\Prod\to\GL(2,\Lambda)$. The Artin
generators $\Gs_1$, $\Gs_2$ corresponding to the chosen geometric basis (the
one used to define $\be_1$, $\be_2$) act \via\
\[*
\Gs_1=\bmatrix-t&1\\0&1\endbmatrix,\qquad
\Gs_2=\bmatrix1&0\\t&-t\endbmatrix,
\]
and an element $\Ga\in\Inn\AA\cong\AA$ maps to $t^{\deg\Ga}\id$.
For the image $\bGb$ of $\Gb\in\Prod$, one has
$\det\bGb=t^{\bdeg\Gb}$.
The restriction $\BG3\to\GL(2,\Lambda)$
is called the
\emph{\rom(reduced\rom) Burau representation}, see~\cite{Burau}; it is
faithful, and for this reason we identify braids and their images in
$\GL(2,\Lambda)$. The image in $\GL(2,\Lambda)$ of the whole group $\Prod$ is
called the \emph{Burau group}~$\Bu3$; it is the central product
\[*
\Bu3=\BG3\odot\<t\id\>=\BG3\times\<t\id\>/\{\Delta^2=t^3\id\}.
\]
The center of $\Bu3$ is the infinite cyclic subgroup generated by the scalar
matrix $t\id$.

Given two submodules $\CU,\CV\subset\AM$, we say that $\CU$ is
\emph{conjugate} to~$\CV$, $\CU\sim\CV$ (\emph{subconjugate} to~$\CV$,
$\CU\prec\CV$) is $\CU=\Gb\CV$ (respectively, $\CU\subset\Gb\CV$) for some
$\Gb\in\BG3$.
In this definition, $\BG3$ can be replaced with the Burau group~$\Bu3$.
Similar terminology is used for subgroups of~$\BG3$ and $\Bu3$: a subgroup
$H\subset\Bu3$ is said to be \emph{conjugate} (\emph{subconjugate}) to a
subgroup $G\subset\Bu3$ if $H=\Gb G\Gb\1$ (respectively,
$H\subset\Gb G\Gb\1$) for some element $\Gb\in\BG3$.

Specializing all matrices at $t=-1$, we obtain epimorphisms
$\BG3\subset\Bu3\onto\tMG:=\SL(2,\Z)$, which factor further to the
\emph{modular representation}
\[*
\pr_\MG\:\BG3,\Bu3\onto\MG:=\PSL(2,\Z)=\tMG/{\pm\id}.
\]
We abbreviate $\bar H:=\pr_\MG H$ and $\bGb:=\pr_\MG\Gb$ for a subgroup
$H\subset\Bu3$ and an element $\Gb\in\Bu3$.
The degree homomorphisms $\bdeg\:\BG3\to\Z$ and $\bdeg\:\Bu3\to\Z$ descend to
epimorphisms $\bdeg\:\MG\onto\CG6$ and ${\bdeg}\bmod2\:\MG\onto\CG2$,
respectively;
the
former coincides with the abelianization epimorphism
$\MG\onto\MG/[\MG,\MG]\cong\CG6$.
The kernels of the projections $\pr_\MG\:\BG3,\Bu3\onto\MG$ are the centers of the
corresponding groups.

\subsection{Skeletons}\label{s.skeletons}
Recall
that a \emph{bipartite graph} is a graph whose vertices are divided into two
kinds, \black-- and \white--, so that the two ends of each edge are of the
opposite kinds. A \emph{ribbon graph} is a graph equipped with a
distinguished cyclic order (\iq. transitive $\Z$-action) on the star of each
vertex. Any graph embedded into an oriented surface~$S$ is a ribbon graph,
with the cyclic order induced from the orientation of~$S$. Conversely, any
finite ribbon graph defines a unique, up to homeomorphism, closed oriented
surface~$S$ into which it is embedded: the star of each vertex is embedded
into a small oriented disk (it is this step where the cyclic order is used),
these disks are connected by oriented ribbons along edges
producing a tubular neighborhood of the graph, and finally each
boundary component of the resulting compact surface is patched with a disk.
(Intuitively, the boundary components patched at the last step are
the \emph{regions} reintroduced
combinatorially in \autoref{def.region} below, where the construction of~$S$
is discussed in more details.)
The surface~$S$ thus constructed is called the \emph{minimal supporting
surface} of the ribbon graph.

Below,
we redefine a certain class of bipartite ribbon
graphs in purely combinatorial terms, relating them to the modular group
and its subgroups.
In spite of this combinatorial approach, we will
freely use the topological language applicable to the geometric
realizations of the graphs.

As is well known, the \emph{modular group} $\MG:=\PSL(2,\Z)$
is generated by two elements
$\X:=(\Bs_2\Bs_1)\1$ and $\Y:=\Bs_2\Bs_1^2$, the defining relations being
$\X^3=\Y^2=1$; thus, $\MG\cong\CG3*\CG2$.
According to~\cite{degt:monodromy} (see also~\cite{Kulkarni}, where this
construction appeared first),
a subgroup $G\subset\MG$ can be described by its
\emph{skeleton} $\Sk:=\Sk_G$, which is the bipartite ribbon graph,
possibly infinite,
defined as
follows: the set of edges of~$\Sk$ is the left $\MG$-set $\MG/G$, its
\black-- and \white-vertices are the orbits of~$\X$ and~$\Y$, respectively,
and the cyclic order (the ribbon graph structure) at each trivalent
\black-vertex is given by the action of~$\X\1$.
The incidence map assigns to an edge~$e$ its \black-- and \white-ends by
sending~$e$ to, respectively, the $\X$- and $\Y$-orbits containing~$e$.
Note that a cyclic order at each mono- or bivalent vertex is unique and hence
redundant; however, it is convenient to agree that the cyclic order at
\emph{all} \black-vertices is given by $\X\1$, and that at \white-vertices is
given by~$\Y$.
The graph $\Sk_G$ is finite if and only if $G$ is a subgroup of finite index.

By definition, $\Sk$ is a connected bipartite ribbon graph,
possibly infinite,
with the following properties:
\roster*
\item
the valency of each \black-vertex is~$3$ or~$1$;
\item
the valency of each \white-vertex is $2$ or~$1$.
\endroster
Such a graph is called an \emph{\rom(abstract\rom) skeleton};
its set of edges can be
regarded as a transitive left $\MG$-set, with the action of $\X\1$ and~$\Y$
given by the cyclic order at the \black-- and \white-vertices, respectively.
A skeleton is \emph{regular} if it has no monovalent vertices.

The skeleton~$\Sk_G$ of a subgroup $G\subset\MG$
is equipped with a distinguished edge~$e$, namely the coset $G/G$.
Conversely, any pair $(\Sk,e)$, where $\Sk$ is a skeleton and
$e$ is a distinguished edge,
gives rise to a subgroup $G:=\stab e\subset\MG$. If no edge is distinguished,
$\Sk$ defines a conjugacy class of subgroups of~$\MG$, which is denoted by
$\Stab\Sk$.

Topologically, we regard a skeleton~$\Sk$ as an orbifold, assigning to each
monovalent \black-- or \white-vertex the ramification index~$3$ or~$2$,
respectively.
Under this convention, the homotopy classes of paths in~$\Sk$ (starting and
ending inside an edge) can be identified with pairs $(e_0,g)$, where
the \emph{initial point}~$e_0$ is an edge of~$\Sk$ and $g\in\MG$; then the
\emph{terminal point} is the edge $e_1:=g\1e_0$.
Hence, we have an isomorphism
\[*
G=\stab e=\piorb(\Sk_G,e),
\]
where the basepoint for the fundamental group is chosen inside the edge~$e$.

\definition\label{def.region}
A \emph{region} in a skeleton~$\Sk$ is an orbit of $\X\Y$. The cardinality of
a region~$R$ is called its \emph{width} and denoted by $\rdeg R$.
A region~$R$ of width~$n$ is also referred to as an \emph{$n$-gon} or
\emph{$n$-gonal region}. (The \black-ends of the edges constituting~$R$ can
be regarded as its corners.)
The region containing an edge~$e$ is denoted by $\reg(e)$.
\enddefinition

Let $\Sk$ be a finite skeleton. Patching each region of~$\Sk$ with an
oriented disk, one obtains
the \emph{minimal supporting surface} $\Supp\Sk$:
it is an oriented closed surface containing~$\Sk$ and
inducing its ribbon graph structure. (More precisely, the
boundary of the disk patching a region $R=\{e_0,e_1,\ldots\}$ is composed by
$e_0,e_0',e_1,e_1',\ldots$, where $e_i':=\Y e_i$; the edges~$e_i'$
appear in the
boundary with the opposite orientation.) The genus~$g$ of $\Supp\Sk$
is called the
\emph{genus} of~$\Sk$. If $\Sk=\Sk_G$ for a finite index subgroup
$G\subset\MG$, then $g$ is also called the \emph{genus} of~$G$; this
definition is equivalent to the conventional one in terms of modular curves,
see~\cite{degt:monodromy}.

The skeleton~$\Sk_G$ and genus of a finite index subgroup $G\subset\Bu3$
(or $G\subset\BG3$)
are
defined as those of the image $\bG\subset\MG$.
Since an inclusion of subgroups
gives rise to
a ramified covering of the
minimal supporting surfaces of their skeletons, one has
\[
\QOPNAME{genus}(H)\ge\QOPNAME{genus}(G)
\quad\text{whenever}\quad
H\prec G.
\label{eq.genus}
\]

\subsection{The type specification}\label{s.type}
Define the \emph{depth} $\depth G$ of a subgroup $G\subset\Bu3$
as the degree of the
minimal positive generator of the cyclic subgroup $G\cap\Ker\pr_\MG$,
or zero if the latter intersection is trivial.
Clearly, $\depth G=0\bmod2$ for any subgroup $G\subset\Bu3$, and
$\depth G=0\bmod6$ if $G\subset\BG3$.

Our primary concern are subgroups of genus zero.
Let $G\subset\Bu3$ be such a subgroup, and denote by $S\scirc:=\Supp\scirc G$
the punctured surface obtained from the sphere $\Supp\Sk_G$ by removing the
center of each region of~$\Sk_G$ and each monovalent vertex of~$\Sk_G$.
Then there is an epimorphism
\[
\pi_1(S\scirc,e)\onto G/\<t^k\id\>,\quad 2k:=\depth G,
\label{eq.pi1.Scirc}
\]
which is included into the commutative diagram
\[*
\CD
\pi_1(S\scirc,e)@>>> G/\<t^k\id\>\\
@VVV@V\cong V\pr_\MG V\\
\piorb(\Sk_G,e)@>\cong>>\bG\rlap.
\endCD
\]
As above, the basepoint for all fundamental groups is chosen inside the
distinguished edge~$e$ of~$\Sk_G$.

Since $S\scirc$ is a punctured sphere, the group $\pi_1(S\scirc,e)$, and
hence also the quotient $G/\<t^k\id\>$, is generated by (the images of)
a system of lassoes in~$S\scirc$ about
the centers of the regions of~$\Sk_G$ and its monovalent vertices.
It follows that the subgroup $G/\<t^k\id\>\subset\Bu3/\<t^k\id\>$
can be described by
means of its \emph{type specification} $\type$, which is a function on the
set of regions and monovalent vertices of $\Sk_G$, taking values in
$\CG{\depth G}$
(with the convention that
$\CG{0}=\Z$) and defined as follows: the value of $\type$ on a
monovalent vertex or a region is the degree of the lift to $G/\<t^k\id\>$ of
the corresponding lasso about the vertex or the center of the region,
respectively. This function is well defined and has the following properties.

\proposition[see~\cite{degt:Alexander}]\label{type}
Let $d=6$ if $G\subset\BG3$ and $d=2$ otherwise. Then\rom:
\roster
\item\label{tp.dp}
$\depth G=0\bmod d$\rom;
\item\label{tp.R}
$\type(R)=\rdeg R\bmod d$ for any region~$R$\rom;
\item\label{tp.black}
$\type(\BLACK)=2\bmod d$ and $3\type(\BLACK)=0$\rom;
\item\label{tp.white}
$\type(\WHITE)=3\bmod d$ and $2\type(\WHITE)=0$\rom;
\item\label{tp.sum}
the sum of all values of~$\type$ equals zero.
\endroster
Given a skeleton~$\Sk$,
a pair $(\depth,\type)$ satisfying
conditions~\iref{tp.dp}--\iref{tp.sum}
above
defines a unique subgroup $G\subset\Bu3$\rom;
one has $G\subset\BG3$ if and only if the pair $(\depth,\type)$
satisfies conditions~\iref{tp.dp}--\iref{tp.white} with $d=6$.
\pni
\endproposition

\section{Trigonal curves}\label{S.tc}

\subsection{Trigonal curves in Hirzebruch surfaces}\label{s.tc}
A \emph{Hirzebruch surface} is a geometrically ruled rational surface
$\pp\:\Sigma_d\to\BB\cong\Cp1$ with an exceptional section~$\EE$ of
self-intersection $-d\le0$. If $d>0$, such a section is unique.
A \emph{\rom(generalized\rom) trigonal curve} is a reduced curve
$\CC\subset\Sigma_d$,
not containing~$\EE$ or a fiber of~$\pp$ as a component, and
such that the restriction $\pp\:\CC\to\BB$ is a map of degree three.
A trigonal curve is \emph{genuine} or \emph{proper} if it is disjoint from
the exceptional section~$\EE$.
A \emph{singular fiber} of a trigonal curve~$\CC$ is a fiber of~$\pp$
intersecting $\CC\cup\EE$ at fewer that four points.

A \emph{positive \rom(negative\rom) Nagata transformation} is a birational
map $\Sigma_d\dashrightarrow\Sigma_{d\pm1}$ consisting in blowing up a
point~$P$ in (respectively, not in)
the exceptional section $\EE$ and blowing down the proper transform
of the fiber through~$P$.
A \emph{$d$-fold Nagata transformation} is a sequence of $d$ Nagata
transformations in the same fiber and \emph{of the same sign}. Two trigonal
curves are \emph{Nagata equivalent} (\emph{$d$-Nagata equivalent}) if they
can be related by a sequence of Nagata transformations (respectively,
$d$-fold Nagata transformations).

By a sequence of positive Nagata transformations, any trigonal curve~$\CC$
can be made proper; the result is called a \emph{proper model} of~$\CC$.

In appropriate affine coordinates $(x,y)$ in~$\Sigma_d$ such that
$E=\{y=\infty\}$, a proper trigonal curve~$C$ can be given by its
\emph{Weierstra{\ss} equation}
\[
y^3+g_2(x)y+g_3(x)=0,
\label{eq.Weierstrass}
\]
where $g_2,g_3$ are certain polynomials in~$x$. The
\emph{\rom(functional\rom) $j$-invariant} of~$C$ is the meromorphic function
$j_C\:B\to\Cp1=\C\cup\{\infty\}$ given by
\[*
j_C(x)=-\frac{4g_2^3}\Delta,\qquad
\text{where}\quad \Delta:=-4g_2^3-27g_3^2
\]
is the discriminant of~\eqref{eq.Weierstrass} with respect to~$y$.
(We use Kodaira's normalization, with respect to which the `special' values
of the $j$-invariant are $0$, $1$, and~$\infty$.)
By definition, $j_C$ is preserved by Nagata transformations, and the
$j$-invariant of an improper trigonal curve is defined as that of any of its
proper models. A
curve~$C$ is called \emph{isotrivial} if
$j_C=\const$.

\subsection{The monodromy group}\label{s.monodromy}
In this subsection, we outline the construction
and basic properties of the braid monodromy of a
trigonal curve. For more details and all proofs, which are omitted here,
we refer to~\cite{degt:e6} and~\cite{degt:dessin}.

Let $C\subset\Sigma_d\to B$ be a proper trigonal curve.
A \emph{monodromy domain}
is a closed topological disk
$\disk\subset B$ containing in its interior all singular
fibers of~$C$. A continuous section $s\:\disk\to\Sigma_d$ of~$\pp$ is called
\emph{proper} if its image is disjoint from both~$E$ and the fiberwise convex
hull of~$C$ (with respect to the canonical affine structure in the
\emph{affine fibers} $F_b^\circ:=\pp\1(b)\sminus E$, $b\in B$,
which are affine spaces over~$\C$).
Since $\disk$ is contractible, a proper section exists and is unique up to
homotopy in the class of such sections.

Fix a monodromy domain~$\disk$ and a proper section~$s$ over~$\disk$. Let
$b_1,\ldots,b_r\in\disk$ be the singular fibers of~$C$, and denote
$\disk\scirc:=\disk\sminus\{b_1,\ldots,b_r\}$. Then, $s$ is a section of the
locally trivial fibration
$\pp\:\pp\1(\disk\scirc)\sminus(C\cup E)\to\disk\scirc$, and the monodromy of
the associated bundle with the discrete fibers
$\Aut\pi_1(F_b^\circ\sminus C,s(b))$, $b\in\disk\scirc$, gives rise to an
anti-homomorphism
$\bm\:\pi_1(\disk\scirc,b)\to\Aut\pi_F$, where
$\pi_F:=\pi_1(F_b^\circ\sminus C,s(b))$, $b\in\disk\scirc$, is the
fundamental group of a fixed nonsingular affine fiber punctured at~$C$.
The latter anti-homomorphism
is called the \emph{braid monodromy} of~$C$,
and its image $\BM_C:=\Im\bm\subset\Aut\pi_F$ is
called the \emph{monodromy group} of~$C$.

The free group $\pi_F$ has a distinguished class of geometric bases; a choice
of one of these bases identifies $\pi_F$ with $\AA$. (In fact, if
$j_C(b)\ne0,1$, then $\pi_F$ has a \emph{canonical basis}
$\{\Ga_1,\Ga_2,\Ga_3\}$, which is well defined up to conjugation by
$\Gr:=\Ga_1\Ga_2\Ga_3$.) Under this identification, the monodromy~$\bm$
takes values in the braid group $\BG3\subset\Aut\AA$ and, up to conjugation
in~$\BG3$, the monodromy group $\BM_C$ is independent of the choices made in
the construction.

The following statement is crucial for \autoref{th.main}.

\theorem[see \cite{degt:trigonal}]\label{th.monodromy}
The monodromy group of a non-isotrivial proper trigonal curve is of genus
zero. Conversely, given a subgroup $G\subset\BG3$ of genus zero
and depth $6d>0$, there is a unique, up to isomorphism and
$d$-Nagata equivalence, proper trigonal curve $\CC_G$ such that,
for another non-isotrivial proper trigonal curve~$\CC$, one has
$\BM_\CC\prec G$ if and only if $\CC$ is $d$-Nagata equivalent to a
curve induced from~$\CC_G$.
This curve $\CC_G$
is called the
\emph{universal curve} corresponding to~$G$.
\pni
\endtheorem

Now, let $C$ be an improper trigonal curve. Consider a proper model~$C'$
of~$C$ and, after making the necessary choices, its braid monodromy
$\bm'\:\pi_1(\disk\scirc,b)\to\Aut\pi_F$.
Let $\{\Gg_1,\ldots,\Gg_r\}$ be a geometric basis for the free group
$\pi_1(\disk\scirc,b)$. To each basis element~$\Gg_j$ one can assign
the \emph{slope} $\slope_j\in\pi_F$, which depends on both curves~$C$, $C'$
and the generator~$\Gg_j$.
In this notation, the \emph{braid monodromy} of~$C$ is defined as
the anti-homomorphism $\bm\:\Gg_j\mapsto\bm_j$, $j=1,\ldots,r$,
where $\bm_j$ is the
automorphism $\Ga\mapsto\slope_j\1\bm(\Gg_j)\slope_j$, $\Ga\in\pi_F$.
The image $\BM_C:=\Im\bm$ is called the \emph{monodromy group} of~$C$; under
the identification $\pi_F=\AA$ it is a subgroup of $\Prod$.

\subsection{The Zariski--van Kampen theorem}\label{s.vanKampen}
The following theorem is the most well-known means of computing the
fundamental group of the complement of an algebraic curve. It is essentially
contained in~\cite{vanKampen}. There is a great deal of modifications and
generalizations of this theorem making use of various pencils;
the particular case of improper trigonal
curves is treated in details in~\cite{degt:dessin}.

\theorem[see~\cite{degt:dessin}]\label{th.vanKampen}
Let $C\subset\Sigma_d$ be a trigonal curve, and let
$\BM_C\subset\Prod$ be its monodromy group. Then one has a presentation
\[*
\piaff{C}=\AA/\{\Ga=\Gb(\Ga),\ \Ga\in\AA,\ \Gb\in\BM_C\}.
\def\qedsymbol{\pnisymbol}\pushQED{\qed}\qedhere
\]
\endtheorem

A presentation of the group $\piaff{C}$ as in \autoref{th.vanKampen} is
called \emph{geometric}.

\subsection{The Alexander modules}\label{s.Alexander}
Given a subgroup $G\subset\Prod$, let
\[*
\bVM_G:=\sum_{\Gb\in G}\Im(\Gb-\id)\subset\AM,\qquad
\VM_G:=\sum_{\Gb\in G,\ \Ga\in\AA}\Lambda[\Gb(\Ga)\cdot\Ga\1]\subset\AM
\]
and define the \emph{Alexander module} $\AM_G:=\AM/\VM_G$ and the
\emph{extended Alexander module} $\bAM_G:=\AM/\bVM_G$.
As in the case of curves,
pick an algebraic number $\xi\in\bk$ over~$\Bbbk_p$,
consider the specializations
\[*
\bAM_G(\xi):=(\bAM_G\otimes\Bbbk_p)/\mp_\xi,\qquad
\AM_G(\xi):=(\AM_G\otimes\Bbbk_p)/\mp_\xi,
\]
and define the subspaces
\[*
\bVM_G(\xi):=\Ker[\AM(\xi)\onto\bAM_G(\xi)],\qquad
\VM_G(\xi):=\Ker[\AM(\xi)\onto\AM_G(\xi)].
\]
Clearly,
\[*
\bVM_G(\xi)=\sum_{\Gb\in G}\Im(\Gb(\xi)-\id)\subset\AM(\xi),
\]
where $\Gb\mapsto\Gb(\xi)$ is the composition of the Burau representation and
specialization homomorphism
$\GL(2,\Lambda)\to\GL(2,\bk)$.
In particular,
both $\bVM_G$ and $\bVM_G(\xi)$ depend on the image of $G$ in $\Bu3$ only and
thus can be defined for subgroups of~$\Bu3$.

\lemma[see~\cite{degt:Alexander}]\label{lem.bVM=VM}
For any subgroup $G\subset\Prod$ and any algebraic number~$\xi$,
one has $\bVM_G(\xi)\subset\VM_G(\xi)$\rom; hence, there is an epimorphism
$\bAM_G(\xi)\onto\AM_G(\xi)$.
If $G\subset\BG3$ and $\xi^2+\xi+1\ne0$, then
$\bVM_G(\xi)=\VM_G(\xi)$ and $\bAM_G(\xi)=\AM_G(\xi)$.
\pni
\endlemma

According to \autoref{th.vanKampen}, for a trigonal curve~$C$
and algebraic number~$\xi$
one has
$\AM_C(\xi)=\AM_G(\xi)$, where $G:=\BM_C$; the corresponding epimorphism
$\AM(\xi)\onto\AM_C(\xi)$ is called a \emph{geometric presentation} of the
Alexander module of~$C$.
Hence, there is an epimorphism
$\bAM_G(\xi)\onto\AM_C(\xi)$, and \autoref{th.main} is essentially a
consequence
of the following restatement in terms of the monodromy groups.

\theorem\label{th.groups}
Let $G\subset\Bu3$ be a subgroup of genus zero and let
$\xi\in\bk\supset\Bbbk_p$ be an algebraic number such that $\bAM_G(\xi)\ne0$.
Then $N:=\ord(-\xi)<\infty$. Furthermore, one has $N\le6$ unless
$(p,\mp_\xi)$ is one of the pairs listed in \autoref{tab.factors}.
Each pair listed in the table is realized by a certain subgroup
$G\subset\Bu3$ of genus zero\rom; the pairs marked with a $^*$ are also
realized by subgroups $G\subset\BG3$ of genus zero.
\endtheorem

This theorem is proved in \autoref{S.proof.main}, see \autoref{s.groups}.

\subsection{The universal subgroups}\label{s.universal}
The existence part of \autoref{th.groups} is based on the concept of
universal subgroup. Fix an algebraic number~$\xi$ and consider a subspace
$\VM\subset\AM(\xi)$. Then the subset
\[*
G_\VM:=\{\Gb\in\Bu3\,|\,\Im(\Gb(\xi)-\id)\subset\VM\}
\]
is a subgroup of~$\Bu3$; it is called the \emph{universal subgroup}
corresponding to~$\VM$.

\remark\label{rem.universal}
Clearly, one has $G_0=\Ker[\Gb\mapsto\Gb(\xi)]$ and
$G_\AM=\Bu3$. In all other cases, $\VM=\bk\bv$
for a certain vector $\bv=a_1\be_1+a_2\be_2\in\AM(\xi)$ and the universal
subgroup $G_\VM$ is given by linear equations:
$\Gb\in G_\VM$ if and only if $\bv^\perp\Gb(\xi)=\bv^\perp$, where
$\bv^\perp:=[a_2,-a_1]$ generates the annihilator
$\VM^\perp\subset\AM(\xi)^*$.
\endremark

The following statements are obvious:
\roster
\item\label{ug.1}
if $G:=G_\VM$, then $\bVM_G\subset\VM$;
\item\label{ug.2}
one has $\bVM_G\prec\VM$ if and only if $G\prec G_\VM$.
\endroster
Here, in Statement~\iref{ug.1}, the inclusion may be proper; in fact, very few subspaces
of dimension one result in nontrivial universal subgroups, \cf.
\autoref{cor.types}.

\lemma\label{lem.width|N}
Let $N:=\ord(-\xi)$, and assume that $2\le N<\infty$
and that $G:=G_\VM$ is the universal subgroup corresponding to a subspace
$\VM\subset\AM(\xi)$.
Then,
the width of each region of the skeletons
$\Sk_G$ and $\Sk_{G\cap\BG3}$ divides~$N$.
\endlemma

\proof
Observe that
\[*
\Gs_1^N=\bmatrix(-t)^N&\tilde\Gf_N(-t)\\0&1\endbmatrix,
\]
where $\tilde\Gf_N(t):=(t^N-1)/(t-1)$. Hence, $\Gs_1^N(\xi)=\id$
and $(\X\Y)^N\in G_\VM\cap\BG3$.
\endproof

\subsection{Digression: the case \pdfstr{N<=6}{$N\le6$}}\label{s.N<=6}
For completeness, we discuss a few extensions of the results
of~\cite{degt:Alexander} concerning the specializations of the Alexander
modules at algebraic numbers~$\xi$ with $N:=\ord(-\xi)\le6$.

Strictly speaking, only irreducible
curves
(equivalently, subgroups of $\Bu3$ with transitive image in~$\SG3$)
are considered in~\cite{degt:Alexander}. However,
the preliminary results of~\cite{degt:Alexander} hold in the general case.
Thus, if $2\le N\le5$ and $G_\VM$ is the universal subgroup corresponding to a
submodule
$\VM\subset\AM/\cp_N(-t)$ (where $\cp_N$ is the
cyclotomic polynomial of order~$N$), then $\bG_\VM\subset\MG$ is a congruence
subgroup of level~$N$.
(In fact, this statement is contained in \autoref{lem.width|N}, as the
principal congruence subgroup of level $N\le5$ is of genus zero and is
normally generated by $\Bs_1^N$.)
The number of such subgroups is finite and, using, \eg, the tables found
in~\cite{Cummins.Pauli} and trying various type specifications, one arrives
at a finite list of submodules of the form $\bVM_G\subset\AM/\cp_N(-t)$.
Details are left to the reader, and the final result, in terms of the
specializations $\bAM_G(\xi)$, is represented in \autoref{tab.factors.<5}.
Listed in the table are:
\roster*
\item
the values of $p$, $N$, and $\mp_\xi$,
\item
the corresponding subspace
$\bVM_G(\xi)\subset\AM(\xi)$ (see below),
\item
the projection $\bG\subset\MG$ of the corresponding universal
subgroup~$G$,
in the notation of~\cite{Cummins.Pauli} and, whenever available, in the
conventional notation,
and
\item
a list of dependencies, \ie, whether the non-vanishing of the module
$\AM_G(\xi)$ implies the non-vanishing of another module $\AM_G(\xi')$ for
the same group~$G$.
\endroster
The subspace
$\bVM_G(\xi)\subset\AM(\xi)$ is
either~$0$ or conjugate to $\bk\bv_T$, where $T$ is the \emph{type} $\vI$,
$\vII$, $\vIII$, or~$\vIV$ listed in the table and
$\bv_T:=a_T(\xi)\be_1+\be_2$, see~\eqref{eq.types}.
The implications in the last column are given by the inclusions of the
universal subgroups, see~\cite{Cummins.Pauli}.

\table[t]
\caption{Alexander modules $\AM_C(\xi)$ with $N:=\ord(-\xi)\le5$%
 }\label{tab.factors.<5}
\def\fref#1{\ref{#1}}
\def\>{$\Rightarrow\;$}
\def\={$\Leftrightarrow\;$}
\setcounter{line}0
\let\tab\quad
\def\linelabel#1{\refstepcounter{line}\theline\label{#1}}
\centerline{\vbox{\halign{\tab\hss\linelabel{#}\tab&
 \tab\hss$#$\hss\tab&\tab\hss$#$\tab&\tab$#$\hss\tab&\tab\hss$#$\hss\tab&
 \tab$#$\hss\tab&\tab#\hss\tab&\tab#\hss\tab\cr
\omit\tab\hss\#\tab&p&$N$&\mp_\xi\in\FF{p}[t]&\bVM_G&
 \bG\subset\MG&Remarks\cr
\noalign{\vskip2pt \hrule\vskip3pt}
0,2&0&2&t-1&\vI&2B^0=\MG_1(2)&
  \=\ref{2,1}\cr
 0,2,0&&&&0&2C^0=\MG(2)&
  \=\ref{2,1,0}\cr
 0,3&&3&t^2-t+1&\vI&3B^0=\MG_1(3)&
  \=\ref{3,1}\cr
 0,3,0&&&&0&3D^0=\MG(3)&
  \=\ref{3,1,0}\cr
 0,4&&4&t^2+1&\vI&4B^0=\MG_1(4)&
  \>\ref{0,2}, \ref{2,1}\cr
 &&&&0&4G^0=\MG(4)&
  \>\ref{0,2,0}, \ref{2,1,0}\cr
 0,5&&5&\cp_5(-t)&\vI&5D^0=\MG_1(5)&
  \>\ref{5,1}\cr
 0,5,0&&&&0&5H^0=\MG(5)&
  \>\ref{5,1,0}\cr
2,1&2&1&t+1&\vI&2B^0=\MG_1(2)&
  \=\ref{0,2}\cr
 2,1,0&&&&0&2C^0=\MG(2)&
  \=\ref{0,2,0}\cr
 2,3&&3&t^2-t+1&\vII&3A^0=\MG^3\cr
 2,5&&5&\cp_5(-t)&\vIV&5E^0\cr
3,1&3&1&t+1&\vI&3B^0=\MG_1(3)&
  \=\ref{0,3}\cr
 3,1,0&&&&0&3D^0=\MG(3)&
  \=\ref{0,3,0}\cr
 3,2&&2&t-1&\vII&2A^0=\MG^2\cr
 3,4&&4&t^2+1&\vIII&4D^0&
  \>\ref{3,2}\cr
 3,5&&5&\cp_5(-t)&\vIII&5F^0\cr
5,1&5&1&t+1&\vI&5D^0=\MG_1(5)\cr
 5,1,0&&&&0&5H^0=\MG(5)\cr
&7&1&t+1&\vI&7E^0=\MG_1(7)\cr}}}
\small\medskip
\endtable

The case $N=1$ (the maximal dihedral quotients of the fundamental group) is
settled in~\cite{degt:trigonal}:
in this case, the universal subgroups are also congruence subgroups of
$\tMG:=\SL(2,\Z)$ (but not necessarily of level~$1$).

Finally, if $N=6$, the $\BG3$-action on the module $\AM':=\AM/(t^2+t+1)$ has
invariant vector $\bv:=-t\be_1+\be_2$. Hence, in the basis $\{\bv,\be_2\}$,
the $\Bu3$-action is given by upper triangular matrices and can easily be
studied. Assume that $G\subset\Bu3$ is a subgroup of genus zero and the
submodule $\bVM_G\subset\AM'$ is distinct from~$\AM'$. If $G\subset\BG3$
(proper trigonal curves), then $\bVM_G\sim\Lambda'\bu+J\bv$, where
$\Lambda':=\Lambda/(t^2+t+1)$, $\bu$ is one of the following five vectors
\[*
\bu_1:=\be_2,\quad
\bu_2:=(t+2)\be_2,\quad
\bu_3:=2\be_2,\quad
\bu_\circ:=t\be_1+\be_2,\quad
\bu_\bullet:=\be_1-\be_2,
\]
and $J\subset\Lambda'$ is an ideal \emph{of finite index}.
If $G\not\subset\BG3$ (improper curves), then $\bVM_G$ is conjugate to the
submodule generated by one of the following seven (pairs of) vectors:
\[*
\gathered
2\be_2,\bv;\quad
\bv;\quad
(t-1)\be_2,\bv;\quad
2\be_2,(t-1)\bv;\quad
\be_2,(t-1)\bv;\\
(t-1)\be_2,(t-1)\bv;\quad
(t-1)\be_2-\bv,(t-1)\bv.
\endgathered
\]
Details will appear elsewhere.

\section{Proof of Theorems~\ref{th.main} and~\ref{th.groups}}\label{S.proof.main}

\subsection{The set-up}\label{s.set-up}
Fix a subgroup $G\subset\Bu3$
\emph{of genus zero} and let $\Sk:=\Sk_G=\MG/\bG$ be its
skeleton, $e:=\bG/\bG$ the distinguished edge of~$\Sk$,
and $\type$ the type specification of~$G$.
Fix, further, a value~$p$, prime or zero, and an algebraic number
$\xi\in\bk_p$. We assume that $N:=\ord(-\xi)\ge7$; in particular,
$\xi\ne\pm1$ and $\xi^2+\xi+1\ne0$.

We will also make use of the multiplicative order $M:=\ord\xi$.
One obviously has
$M=\ee_p(N)$ and $N=\ee_p(M)$, where $\ee_2(N):=N$ and
\[*
\ee_p(N):=\begin{cases}
2N,&\text{if $N=1\bmod2$},\\
\frac12N,&\text{if $N=2\bmod4$},\\
N,&\text{if $N=0\bmod4$}
\end{cases}
\]
for $p\ne2$ prime or zero.
The $\Bu3$-action on $\AM(\xi)$ factors through $\Bu3/\<t^M\id\>$.
In particular, we can assume that $\depth G=2M$
and pass to the group $G/\<t^M\id\>$.

We are interested in a subgroup~$G$ such that $\bAM_G(\xi)\ne0$. Since
$\dim_{\bk}\AM(\xi)=2$, the latter condition is equivalent to
$\dim_{\bk}\bVM_G(\xi)\le1$
and, according to~\cite{degt:Alexander}, one has $N<\infty$.
A region~$R$
of~$\Sk$ is called \emph{trivial} (\emph{essential}) if $N\divides|\rdeg R$
(respectively, $N\notdivides|\rdeg R$).
Since genus is monotonous, see~\eqref{eq.genus},
we can assume that $G$ is the universal subgroup
corresponding to the subspace $\bVM_G(\xi)\subset\AM(\xi)$.
Then
the width of each region divides~$N$, see \autoref{lem.width|N};
hence, trivial are the regions~$R$
with $\width R=N$, and essential are those with $\width R<N$.

Consider a copy of~$\AA$ and a geometric basis $\Ga_1,\Ga_2,\Ga_3$ with
respect to which the action of~$\BG3$ is given by~\eqref{eq.braid.action}.
Given another edge~$e'$ of~$\Sk$, we fix a path $(e,g)$, $g\in\MG$, from~$e$
to~$e'$, lift~$g$ to an element $\tilde g\in\BG3$, and consider a new
geometric basis
$\Ga_i':=\tilde g(\Ga_i)$, $i=1,2,3$; it is called a \emph{canonical basis}
over~$e'$.
Using these canonical bases for~$\AA$, we define the (local) canonical bases
$\be_1,\be_2$ (over~$e$) and $\be_1',\be_2'$ (over~$e'$) for the universal
Alexander module~$\AM$, see \autoref{s.Burau}.

\subsection{The local modules}\label{s.local}
Consider a region~$R$ or a monovalent vertex~$v$ of~$\Sk$ and denote
$\bVM_*(\xi):=\Im(\bm_*-\id)$, where $\bm_*$ is the monodromy about the
boundary $\partial R$ if $*=R$ or the monodromy about~$v$ if $*=v$.
(More precisely, $\bm_*$ is the image under~\eqref{eq.pi1.Scirc} of a lasso
about the center of~$R$ or~$v$, respectively.)
In view of~\eqref{eq.pi1.Scirc}, one has
$\bVM_G(\xi)=\sum\bVM_*(\xi)$, where $*$ runs over all regions and monovalent
vertices of~$\Sk$.
Hence, a necessary condition for the non-vanishing of $\bAM_G(\xi)$ is
$\dim_{\bk}\bVM_*(\xi)\le1$ for each region and each monovalent vertex.

The submodules $\bVM_*(\xi)$ are easily computed in terms of a local
canonical basis over an edge~$e'$ `close' to the region or vertex in
question.
More precisely, if $*=R$ is a region, we let $e':=\Y e''$, where $e''$ is any
edge contained in~$R$; if $*=v$ is a monovalent \black-vertex, we take
for~$e'$ the only edge incident to~$v$; finally, if $*=v$ is a monovalent
\white-vertex, we let $e'=\X e''$, where $e''$ is the only edge incident
to~$v$.

The following two
statements are contained in~\cite{degt:Alexander}.

\lemma[see~\cite{degt:Alexander}]\label{lem.special}
In the notation above, assume that $\dim_{\bk}\bVM_*(\xi)\le1$, where $*$ is
a region~$R$ or a monovalent vertex~$v$. Let $M:=\ord\xi=\ee_p(N)$.
\roster
\item\label{sp.trivial}
If $R$ is a trivial region, then $\type(R)=\width R\bmod2M$ and
$\bVM_R(\xi)=0$.
\item\label{sp.essential}
Essential regions are subdivided into two types, $\vI$ and~$\vII$,
as explained below.
\item\label{sp.I}
If $R$ is a region of type~$\vI$, then $\type(R)=\width R\bmod2M$ and
$\bVM_R(\xi)=\bk\be_2'$.
\item\label{sp.II}
If $R$ is a region of type~$\vII$ and $n:=\width R$, one has\rom:
if $n$ is even or $p=2$, then $\type(R)=-n\bmod2M$\rom;
otherwise, $\type(R)=M-n\bmod2M$ and $M$ is even\rom;
in both cases, $\bVM_R(\xi)=\bk(\xi\1(\xi+1)\be_1'+\be_2')$.
\item\label{sp.III}
If $v$ is a monovalent \black-vertex, one has\rom:
if $p\ne3$, then $M=0\bmod3$ and $\type(v)=\pm\frac23M\bmod2M$\rom;
otherwise, $M\ne0\bmod3$ and $\type(v)=0\bmod2M$\rom;
in both cases, $\bVM_v(\xi)=\bk(-\xi^s\be_1'+\be_2')$, where
$s:=\frac12\type(v)-1$.
\item\label{sp.IV}
If $v$ is a monovalent \white-vertex, then $M$ is odd,
$\type(v)=M\bmod2M$, and $\bVM_v(\xi)=\bk(\xi^s\be_1'+\be_2')$, where
$s:=\frac12(M-1)$.
\pni
\endroster
\endlemma

\lemma[see~\cite{degt:Alexander}]\label{lem.distinct}
Assume, in addition, that $\bAM_G(\xi)\ne0$. Then\rom:
\roster
\item\label{distinct.3}
at most one of the three regions
incident to a trivalent \black-vertex is essential\rom;
\item
the region incident to a monovalent vertex is trivial\rom;
\item
two monovalent vertices cannot be incident to a common edge.
\pni
\endroster
\endlemma

Note that in
\autoref{lem.distinct}\iref{distinct.3} we do \emph{not}
assume that the three regions are pairwise distinct.
In particular, it follows that a \black-vertex may appear at most once as a
corner of an essential region.

\subsection{The case \pdfstr{V\sb G(xi)<>0}{$\bVM_G(\xi)=0$}}\label{s.VM=0}
If $\bVM_G(\xi)=0$, then, due to \autoref{lem.special}, $\Sk$ is a regular
skeleton and the widths of all regions of~$\Sk$ are multiples of~$N$. Hence,
by Euler's formula
(see, \eg, \eqref{eq.Euler} below), one has $N\le5$. This case is settled
in~\cite{degt:Alexander}, where it is shown that $\bG=\MG(N)$ is the principal
congruence subgroup of level~$N$.

\subsection{The case \pdfstr{0<>V\sb G(xi)<>A(xi)}{$0\subsetneq\bVM_G(\xi)\subsetneq\AM(\xi)$}}
From now on, we assume that
$\dim_{\bk}\bVM_G(\xi)=1$,
\ie, the skeleton~$\Sk$ has at least one essential region or monovalent
vertex.

Consider an edge~$e'$ of~$\Sk$. If $e'$ is the support of a canonical basis
used in the computation of a local module $\bVM_*(\xi)$, see the explanation
prior to \autoref{lem.special}, we assign to~$e'$ a \emph{type}
$T(e')$ as follows:
\roster*
\item
if $*=R$ is an essential region of type~$\vI$ or~$\vII$,
see \autoref{lem.special}\iref{sp.I} and~\iref{sp.II},
then $e'$ is of type~$\vI$ or~$\vII$, respectively;
\item
if $*=v$ is a monovalent \black-vertex and $p\ne3$, then $e'$ is of type
$\vIII_\pm$, where $\type(v)=\pm\frac23M\bmod2M$,
see \autoref{lem.special}\iref{sp.III};
\item
if $*=v$ is a monovalent \black-vertex and $p=3$, then $e'$ is of
type~$\vIII$;
\item
if $*=v$ is a monovalent \white-vertex, then $e'$ is of type~$\vIV$,
see \autoref{lem.special}\iref{sp.IV};
\item
otherwise ($e'$ is not related to a `special' fragment of~$\Sk$), $e'$ is of
type~$0$.
\endroster
An edge of type $T\ne0$ is called \emph{special}.
According to Lemmas~\ref{lem.special} and~\ref{lem.distinct}, the type is
well defined, \ie, an edge cannot be related to two distinct `special'
fragments. (Indeed, otherwise the subspace $\bVM_G(\xi)$ would contain
a pair of
linearly independent vectors and one would have $\bAM_G(\xi)=0$.)
In other words, there is a well defined \emph{surjective} map
\[
\psi\:\Espec\onto\{\text{monovalent vertices}\}\cup
 \{\text{essential regions}\},
\label{eq.special}
\]
where $\Espec$ is the set of the special edges of~$\Sk$.
It follows also that to each special edge $e'$ one can assign the
\emph{local subspace}
$\bVM_{e'}(\xi):=\bVM_{\psi(e')}\subset\AM$. If $e'=\Bb\1e$, $\Gb\in\BG3$,
\autoref{lem.special} implies that
$\bVM_{e'}(\xi)=\bk\bigl(\Gb\bv_{T(e')}\bigr)$,
where $\bv_T:=a_T(\xi)\be_1+\be_2$ and the
Laurent
polynomial $a_T(t)$, $T\ne0$,
is given by
\[
a_{\vI}=0,\quad
a_{\vII}=t\1(t+1),\quad
a_{\vIII}=-t^s,\quad
a_{\vIV}=t^{(M-1)/2}.
\label{eq.types}
\]
Here, $s:=\pm(M/3)-1$ for $T=\vIII_\pm$ ($p\ne3$) and $s:=-1$ for $T=\vIII$
($p=3$).

\corollary[see~\cite{degt:Alexander}]\label{cor.types}
If $G\subset\Bu3$ is a subgroup of genus zero,
$N:=\ord(-\xi)\ge7$,
and $\bVM_G(\xi)\subset\AM$ is
a subspace of dimension one, then one has
$\bVM_G(\xi)\sim\bk\bv_T$ for some type $T\ne0$
\rom(in fact, for any type $T\ne0$ present in the skeleton $\Sk_G$\rom).
\done
\endcorollary

(According to~\cite{degt:Alexander}, the conclusion of \autoref{cor.types}
also holds for $N\le5$; the only exception is the case $N=6$, \ie,
$\xi^2+\xi+1=0$.)

Let $\Rtriv$ be the
set of the \emph{trivial} regions of~$\Sk$.
Let, further, $k_N:=\lceil5/(N-6)\rceil$,
\ie, $k_7=5$, $k_8=3$, $k_9=k_{10}=2$, and $k_N=1$ for $N\ge11$.

\lemma\label{lem.kN}
Assume that there is a map $\Gf\:\Espec\to2^{\Rtriv}$ with the following
properties\rom:
\roster*
\item
$\ls|\Gf(e')|\ge k_N$ for each special edge $e'\in\Espec$\rom;
\item
$\Gf(e')\cap\Gf(e'')=\varnothing$ whenever $e'\ne e''$.
\endroster
Then $G$ is \emph{not} a subgroup of genus zero.
\endlemma

\proof
Let $\nblack$ and $\nwhite$ be the numbers of \emph{monovalent} \black-- and
\white-vertices of~$\Sk$,
and let $n_i$ be the number of its
regions of width $i\ge1$.
As a simple consequence of Euler's formula, $\Sk$ is of genus zero if and
only if
\[
3\nwhite+4\nblack+\sum_{i=1}^N(6-i)n_i=12.
\label{eq.Euler}
\]
Recall that a region~$R$ is trivial if and only if $\width R=N$, \ie,
$\ls|\Rtriv|=n_N$.
Replacing
in~\eqref{eq.Euler}
all coefficients except $(6-N)$ with their
maximum $5=\max_{i\ge1}\{3,4,6-i\}$,
in view of~\eqref{eq.special} we obtain the inequality
$5\ls|\Espec|>(N-6)\ls|\Rtriv|$.
On the other hand, under the hypotheses of the lemma, we have
$\ls|\Rtriv|\ge k_N\ls|\Espec|\ge5\ls|\Espec|/(N-6)$.
\endproof

\subsection{Reduction to a finite number of cases}\label{s.reduction}
We still assume that $G\subset\Bu3$ is the universal subgroup corresponding
to a subspace $\bVM_G(\xi)\subset\AM(\xi)$ of dimension~$1$.

In order to construct a
`universal'
map~$\Gf$ as in \autoref{lem.kN},
we fix a value of~$N$ and
consider a finite set $B=\{\Gb_1,\ldots,\Gb_k\}\subset\Bu3$, $k\ge k_N$,
with all projections $\Bb_i\in\MG$ pairwise distinct.
For a type $T\ne0$, denote $\bv_T(t):=a_T(t)\be_1+\be_2\in\AM$, so that
$\bv_T=\bv_T(\xi)$, and consider the Laurent polynomials
\[*
D_{ij,l}(T',T'')(t):=
 \det\bigl[\Gs_1^l\Gb_i\bv_{T'}(t)\bigm|\Gb_j\bv_{T''}(t)\bigr]\in\Lambda,
\]
where
\[
\gathered
T',T''\ne0,\quad
 i,j=1,\ldots,k,\quad
 l=0,\ldots,N-1,\\
\text{and}\quad T'\ne T''\quad
 \text{or}\quad i\ne j\quad
 \text{or}\quad l\ne0.
\endgathered
\label{eq.indices}
\]
Note that excluded in~\eqref{eq.indices}
are precisely those sequences
$(T',T'',i,j,l)$ for which the determinant is identically zero.

\lemma\label{lem.det}
Let
$e',e''\in\Espec$ be two special edges, not necessarily distinct,
and let $\Gb',\Gb''\in\Bu3$.
Then, if $\Bb'e'=\Bb''e''$, one must have
$\det\bigl[\Gb'\bv_{T(e')}\bigm|\Gb''\bv_{T(e'')}\bigr]=0$.
\endlemma

\proof
Replacing~$G$ with a conjugate subgroup, we can assume that
$\Bb'e'=\Bb''e''$ is the distinguished edge~$e$. Then the vectors
$\Gb'\bv_{T(e')}$ and $\Gb''\bv_{T(e'')}$ span $\bVM_{e'}(\xi)$ and
$\bVM_{e''}(\xi)$, respectively, and,
unless these vectors are linearly dependent, we
have $\dim_{\bk}\bVM_G(\xi)\ge2$, \ie, $\bVM_G(\xi)=\AM(\xi)$.
\endproof

\lemma\label{lem.vanishing}
Assume that $G$ is a subgroup of genus zero and
that $N\ge7$. Then, for any subset
$B\subset\Bu3$ of size $k\ge k_N$, there is a sequence
$(T',T'',i,j,l)$ as in~\eqref{eq.indices} such that
$D_{ij,l}(T',T'')(\xi)=0$.
Furthermore, for at least one of such sequences one has
$\bVM_G(\xi)\sim\bk\bv_{T'}\sim\bk\bv_{T''}$.
\endlemma

\proof
Assume that the conclusion
does \emph{not} hold, \ie, that all determinants are
non-zero. Then, by \autoref{lem.det}, for any pair of special edges
$e',e''\in\Espec$ one has $\Bs_1^l\Bb_ie'\ne\Bb_je''$ whenever $e'\ne e''$ or
$i\ne j$ or $l\ne0\bmod N$. In particular
(from the special case $e'=e''$ and $i=j$), each
region $\reg(\Bb_ie')$ is trivial and, letting
\[*
\Gf(e'):=\bigl\{\reg({\Bb_ie'})\bigm|i=1,\ldots,k\bigr\},
\]
we obtain a well defined map $\Gf\:\Espec\to2^{\Rtriv}$ satisfying the
hypotheses of \autoref{lem.kN}. Hence, $G$ is not of genus zero.

For the last statement, observe that, if $D_{ij,l}(T',T'')(\xi)\ne0$ for all
types $T',T''$ present in $\Sk_G$, then the map~$\Gf$ in this particular
skeleton
is still well defined and satisfies the hypotheses of
\autoref{lem.kN}. Hence, again, $G$ is not of genus zero.
\endproof

Fix a value $N\ge7$, consider a subset
$B:=\{\Gb_1,\ldots,\Gb_k\}\subset\Bu3$, and compute the resultants
$\CR_{ij,l}(T',T'')\in\Z$ of
the determinants $D_{ij,l}(T',T'')(t)$
and the cyclotomic polynomial
$\cp_N(-t)$, where $(T',T'',i,j,l)$ is
an index sequence as in~\eqref{eq.indices}.
The set~$B$ is called \emph{informative} if $k\ge k_N$
and all $\CR_{ij,l}(T',T'')\ne0$ in~$\Z$.
Due to \autoref{lem.vanishing}, the existence of an informative set, see
below, rules out the case $p=0$.
(In~\cite{degt:Alexander}, this case was
prohibited for irreducible curves only.)
Furthermore, each informative set~$B$ gives rise to a finite collection
$\CE(B)$ of `exceptional' triples $(p,\mp_\xi,T)$ such that there \emph{may}
exist a subgroup $G\subset\Bu3$ of genus zero with
$\bVM_G(\xi)\sim\bk\bv_T(\xi)$. This list is obtained as follows: for each
resultant $\CR_{ij,l}(T',T'')\ne\pm1$, we let $T=T'$ and record all prime
divisors~$p$ of $\CR_{ij,l}(T',T'')$
(so that $\CR_{ij,l}(T',T'')=0\bmod p$)
and, for each such divisor~$p$, all irreducible common factors~$\mp_\xi$ of
$D_{ij,l}(T',T'')(t)$ and $\cp_N(-t)$ over~$\Bbbk_p$.

It is shown in~\cite{degt:Alexander} that $N\le26$
and, furthermore, $N\le10$ unless $(p,\mp_\xi)$ is one of the
pairs listed in \autoref{tab.factors}. (Note that the latter
statement can also be
proved using the approach outlined
in this subsection:
for most values $N\ge11$
the subset $B=\{\id\}$ is informative.)
Let $\Gb_1=t\Gs_1\Gs_1\1$ and $\Gb_2=t\Gs_2\1\Gs_1$. (We multiply
the matrices by~$t$ in order to clear the denominators.)
Using \Maple, one can show that each of the following subsets
\begin{alignat*}2
&N=7\:&&\{\id,\Gb_1^2,\Gb_1^3,\Gb_1\Gb_2,\Gb_2\Gb_1\}\
 \text{and}\ \{\id,\Gb_1^2,\Gb_1\Gb_2,(\Gb_1\Gb_2)^2,\Gb_2\Gb_1\};\\
&N=8\:&&\{\id,\Gb_1^2,\Gb_1\Gb_2\}\
 \text{and}\ \{\id,\Gb_1^2,\Gb_1\Gb_2\Gb_1\};\\
&N=9,10\:\quad&&\{\id,\Gb_2\},\ \{\id,\Gb_1\Gb_2\},\
 \text{and}\ \{\id,\Gb_2\Gb_1\}
\end{alignat*}
is informative and, for each subset~$B$, compile the list $\CE(B)$ of
exceptional triples.
(To shorten the further computation, for each~$N$ we consider several
subsets~$B_i$ and take the intersection $\bigcap_i\CE(B_i)$ of the
corresponding lists.) As a result, we obtain a finite list (too long to be
reproduced here) of exceptional triples $(p,\mp_\xi,T)$ that might appear in
the extended Alexander module of a subgroup of genus zero.


\subsection{End of the proof of \autoref{th.groups}}\label{s.groups}
The rest of the proof proceeds as in~\cite{degt:Alexander}: for each
exceptional triple $(p,\mp_\xi,T)$ found in the previous subsection, we
use \Maple\ to
compute the universal subgroup~$G$ of~$\Bu3$ or~$\BG3$ corresponding to the
subspace $\bk\bv_T(\xi)\subset\AM(\xi)$ and select those triples for which
this subgroup is of genus zero.
The result is \autoref{tab.factors}.

For the computation, we specialize the Burau representation at $t=\xi$ and
map $\BG3\subset\Bu3$ to the finite group $\GL(2,\bk)$. (Recall that
$p\ne0$ and $\bk$ is a finite field. In fact, in most cases $\deg\mp_\xi=1$
and hence $\bk=\Bbbk_p$.
In the few exceptional cases, we are working with $(2\times2)$-matrices over
$\Bbbk_p[t]$ considering them modulo~$\mp_\xi$.)
Denote the resulting specialization homomorphism by
$\kk\:\Bu3\to\GL(2,\bk)$.
Then $G\supset\Ker\kk$ and the set of edges of the skeleton $\Sk_G$ is
the quotient of $\Im\kk/\kk(G)$
(or $\kk(\BG3)/\kk(G)$ if the universal subgroup
of~$\BG3$ is to be found)
by the further identification
$m\sim\xi^sm$, where $s\in\Z$ (respectively, $s=0\bmod3$).
The \black-- and \white-vertices of $\Sk_G$ are the orbits of
$\kk(\Gs_2\Gs_1)$ and $\kk(\Gs_2\Gs_1^2)$, respectively, and its regions are
the orbits of $\kk(\Gs_1)$.

Technically, since the image $\Im\kk$ is not known \latin{a priori}, the
coset enumeration proceeds as follows. We start with $m=\id$ and
keep multiplying matrices by $\kk(\Gs_2\Gs_1)$ and $\kk(\Gs_2\Gs_1^2)$,
comparing each matrix against those already recorded.
Each new matrix~$m$ is added to the list together with all products $\xi^sm$,
$s=0,\ldots,M-1$, where $M:=\ee_p(N)$. (If $M=0\bmod3$ and a subgroup
$G\subset\BG3$ is to be found, only the values $s=0\bmod3$ are used.) Note
that the equivalence relation is, in fact,
linear, \cf. \autoref{rem.universal}: for two matrices
$m_1,m_2\in\GL(2,\bk)$ one has $m_1m_2\1\in\kk(G)$ if and only if
$\bv_T^\perp(m_1-m_2)=0$, where $\bv_T^\perp:=[-1, a_T(\xi)]$ generates
the annihilator of the subspace $\bk\bv_T\subset\AM(\xi)$.
This
observation simplifies the coset enumeration.
\qed

\subsection{End of the proof of \autoref{th.main}}\label{s.main}
In view of the epimorphism $\bAM_G(\xi)\onto\AM_C(\xi)$, $G:=\BM_C$, the
restrictions on the pairs $(p,\mp_\xi)$ that may result in a nontrivial
Alexander module follow from
\autoref{th.monodromy} (the monodromy group is a subgroup of genus zero) and
\autoref{th.groups}. If a pair $(p,\mp_\xi)$ can be realized by a subgroup
$G\subset\BG3$ of genus zero
(the lines marked with a $^*$ in \autoref{tab.factors}),
then $\AM_G(\xi)=\bAM_G(\xi)\ne0$, see \autoref{lem.bVM=VM}, and, due to
\autoref{th.monodromy} again, $G$ is the monodromy group of a certain proper
trigonal curve~$C$, so that one has $\AM_C(\xi)=\AM_G(\xi)\ne0$.
\qed

\subsection{Proof of \autoref{ad.<=1}}\label{proof.<=1}
The first statement follows from the computation in \autoref{s.groups}:
in each case resulting in a universal subgroup~$G$
of genus zero, we either start
with a triple $(p,\mp_\xi,T)$ with $T=\vI$ (and hence $\bVM_G=\bk\be_2$) or,
using the coset enumeration, can show that the subspaces
$\bVM_G=\bk\bv_T$ and
$\bk\be_2=\bk\bv_{\vI}$ are conjugate.

The second statement is also proved by a computer aided computation. One
needs to show that, given two universal subgroups $G_1,G_2\subset\Bu3$
corresponding to
two distinct pairs $(p,\mp_\xi)$ and $(q,\mp_\eta)$, the intersection
$G_1\cap G_2'$, where $G_2'\sim G_2$, cannot be of genus zero.
The skeletons $\Sk_i:=\Sk_{G_i}$, $i=1,2$, have already been computed and,
using the double coset formula, one can see that the skeletons of
the intersections of the form $G_1\cap G_2'$, $G_2'\sim G_2$, are the
connected components of the fibered product
$\Sk_1\times_{\SkGamma}\Sk_2$, where $\SkGamma$ is the skeleton of $\MG$
itself.
Considering all products/components one by one, one concludes that they all
have positive genus. Details will appear elsewhere.
\qed

\let\.\DOTaccent
\def\cprime{$'$}
\bibliographystyle{amsplain}
\bibliography{degt}

\end{document}